\documentclass[12pt,oneside,english]{amsart}
\usepackage[latin1]{inputenc}
\usepackage{amssymb}

\makeatletter
\newtheorem{theorem}{Theorem}

\newtheorem{conjecture}{Conjecture}

\newtheorem{corollary}{Corollary}

\usepackage{babel}

\makeatother
\begin{document}
\baselineskip=17pt
\title[Constellations of primes generated by twin primes]{Constellations of primes generated by twin primes}

\author{Vladimir Shevelev}
\address{Department of Mathematics \\Ben-Gurion University of the
 Negev\\Beer-Sheva 84105, Israel. e-mail:shevelev@bgu.ac.il}
 \author{Peter J. C. Moses}
 \address{United Kingdom. e-mail: mows@mopar.freeserve.co.uk}


\begin{abstract}
We study a special set of constellations of primes generated by twin primes.
\end{abstract}

\maketitle

\section{Introduction and main results}
Let $A$ be an infinite sequence of increasing positive integers with the infinite
complementary sequence $\overline{A}.$ Let $H$ be lexicographically first strictly
increasing sequence starting, say, $H(2) = h_0$ with the property that $H(n)\in A$
if and only if $n\in A.$ Shevelev \cite{1} studied several different $H$-sequences
which possess an interesting property: let numbers $a$ and $b$ are arbitrary in $A;$ let
$H_a$ and $H_b$ are defined in the same way as sequence $H$ but with the initial
term $a$ and $b$ respectively; then there exists a position $n=n(a,b)$ in which
$H_a$ and $H_b$ merge. Most likely, it is a difficult problem to find a test for $A,$
when this property holds. Here is only one conjecture in case when $A=P$ is the
sequence of all primes.
 \begin{conjecture}\label{con1} Consider lexicographically first strictly increasing sequence
 starting $H(2) = 3 $ with the property that $H(n)$ is prime if and only if $n$
 is prime. Then for arbitrary primes $a>b\geq3,$  $H_a$ and $H_b$ merge.
 \end{conjecture}
 Note that Conjecture \ref{con1} easily follows from each of the following weaker
 forms: $\alpha)$ $b=3;$ $\beta)$ $a=Nextprime(b), b\geq3.$
 If Conjecture \ref{con1} is true, there appear many natural questions. For example,\newline
 1) For a given primes $a>b\geq3,$ what is the position when $H_a$ and $H_b$ merge?
 \newline
 For the first such pairs $(5,3); (7,3),(7,5); (11,3),(11,5),(11,7);...$  we get the
 following positions:
 \begin{equation}\label{1}
11;\enskip 47,47; \enskip47,47,11;\enskip 47,47,17,17;\enskip683,683,683,683,683;...
\end{equation}
(cf. A276676\cite{2}).\newline
 2) What is the maximal difference between the corresponding terms of $H_a$ and $H_b?$
\newline
 For example, in case when $a$ and $b$ are twin primes: $a=A001359(n-1)+2,\enskip
  b=A001359(n-1), \enskip n\geq2,$ we get
\begin{equation}\label{2}
4, 14, 6, 6, 6, 12, 6, 8, 14, 14, 18, 36, 24, 65, 18, 6, 10, 6, 84, 14, 162, ...
\end{equation}
(cf. A276826\cite{2}).\newline
Let us consider those lesser of twin primes $\{A001359(l)\},\enskip l\geq1,$
which correspond to terms of the sequence (\ref{2}) which not exceeding 6. We have the
sequence $A276848=\{c(n)\}$\cite{2}:
\begin{equation}\label{3}
3, 11, 17, 29, 59, 227, 269, 1277, 1289, 1607, 2129, 2789, 3527, 3917, ...
\end{equation}
We call "nearest" twin primes the corresponding pairs of twin primes $(3,5),$
 $(11,13), (17,19),...$
 A numerical research of sequence (\ref{3}) shows
that we have an interesting phenomenon: although terms c(n) == 7 or 9 (mod 10)
occur often, the first terms c(n)==1 (mod 10) are 11, 165701, ... .
This phenomenon is explained by the following.
\begin{theorem}\label{t1}
$1)$ \enskip $c(n)=p\equiv 7 \pmod {10}$ if and only if, for some $t\geq0,$ $p=30t+17$ and we have the following constellation of five consecutive primes: $p,p+2,p+6,p+12,p+14;$\newline
$2)$ \enskip $c(n)=p\equiv 9 \pmod {10}$ if and only if, for some $t\geq0,$ $p=30t+29$ and we have the following constellation of five consecutive primes: $p,p+2,p+8,p+12,p+14;$\newline
$3)$\enskip$c(n)=p\equiv 1 \pmod {10}$ if and only if, for some $t\geq0,$ $p=30t+11$ and we have the following constellation of seven consecutive primes: $p,p+2,p+6,p+8,p+12,p+18,p+20.$
\end{theorem}
Again in case when $a$ and $b$ are twin primes: $a=A001359(n-1)+2,\enskip
  b=A001359(n-1), \enskip n\geq2,$ consider the smallest $m=m(n)$ such that $H_a(m)-H_b(m)>6$
    or 0 if there is no such $m.$ We get the sequence
\begin{equation}\label{4}
0, 13, 0, 0, 0, 9, 0, 11, 11, 5, 3, 15, 3, 7, 3, 0, 3, 0, 3, 5, 7, 3, 11,...
\end{equation}
\begin{theorem}\label{t2}
The all distinct numbers in the sequence $(\ref{4})$ are $0,3,5,7,9,11,13,15$ and $17$ only.
\end{theorem}

\section{Proof of Theorem 1 (sufficiency)}
1) Let $p=30t+17,$ set $a=30t+19, b=p.$ Let $\{b,a\}$ are twin primes. If also $30t+23, 30t+29, 30t+31$ are prime, then we have\newline
$$H_a(2)=30t+19,\enskip H_b(2)=30t+17;$$
$$H_a(3)=30t+23,\enskip H_b(3)=30t+19;$$
$$H_a(4)=30t+24,\enskip H_b(4)=30t+20;$$
\newpage
$$H_a(5)=30t+29,\enskip H_b(5)=30t+23;$$
$$H_a(6)=30t+30,\enskip H_b(6)=30t+24;$$
$$H_a(7)=30t+31,\enskip H_b(7)=30t+29;$$
$$H_a(8)=30t+32,\enskip H_b(8)=30t+30;$$
$$H_a(9)=30t+33,\enskip H_b(9)=30t+32;$$
$$H_a(10)=30t+34,\enskip H_b(10)=30t+33;$$
$$H_a(11)=Nextprime(30t+34),\enskip H_b(11)=Nextprime(30t+33).$$
Since $Nextprime(30t+34)=Nextprime(30t+33)$ then
$H_a(n)=H_b(n)$ for $n>=11$ and, for $n>=2,$ $\max(H_a(n)-H_b(n)=6$ which
holds for $n=5,6.$\newline

2) Let $p=30t+29,$ set $a=30t+31, b=p.$ Let $\{b,a\}$ are twin primes. If also $30t+37, 30t+41, 30t+43$ are prime,
 then we have\newline
$$H_a(2)=30t+31,\enskip H_b(2)=30t+29;$$
$$H_a(3)=30t+37,\enskip H_b(3)=30t+31;$$
$$H_a(4)=30t+38,\enskip H_b(4)=30t+32;$$
$$H_a(5)=30t+41,\enskip H_b(5)=30t+37;$$
$$H_a(6)=30t+42,\enskip H_b(6)=30t+38;$$
$$H_a(7)=30t+43,\enskip H_b(7)=30t+41;$$
$$H_a(8)=30t+44,\enskip H_b(8)=30t+42;$$
$$H_a(9)=30t+45,\enskip H_b(9)=30t+44;$$
$$H_a(10)=30t+46,\enskip H_b(10)=30t+45;$$
$$H_a(11)=Nextprime(30t+46),\enskip H_b(11)=Nextprime(30t+45).$$
Since $Nextprime(30t+46)=Nextprime(30t+45)$ then
$H_a(n)=H_b(n)$ for $n>=11$ and, for $n>=2,$ $\max(H_a(n)-H_b(n)=6$ which
holds for $n=3,4.$\newline

3) Let $p=30t+11,$ set $a=30t+13, b=p.$ Let $\{b,a\}$ are twin primes. If also $30t+17, 30t+19, 30t+23, 30t+29, 30t+31$ are all prime, then we have\newline
$$H_a(2)=30t+13,\enskip H_b(2)=30t+11;$$
$$H_a(3)=30t+17,\enskip H_b(3)=30t+13;$$
$$H_a(4)=30t+18,\enskip H_b(4)=30t+14;$$
$$H_a(5)=30t+19,\enskip H_b(5)=30t+17;$$
$$H_a(6)=30t+20,\enskip H_b(6)=30t+18;$$
$$H_a(7)=30t+23,\enskip H_b(7)=30t+19;$$
\newpage
$$H_a(8)=30t+24,\enskip H_b(8)=30t+20;$$
$$H_a(9)=30t+25,\enskip H_b(9)=30t+21;$$
$$H_a(10)=30t+26,\enskip H_b(10)=30t+22;$$
$$H_a(11)=30t+29,\enskip H_b(11)=30t+23;$$
$$H_a(12)=30t+30,\enskip H_b(12)=30t+24;$$
$$H_a(13)=30t+31,\enskip H_b(13)=30t+29;$$
$$H_a(14)=30t+32,\enskip H_b(14)=30t+30;$$
$$H_a(15)=30t+33,\enskip H_b(15)=30t+32;$$
$$H_a(16)=30t+34,\enskip H_b(16)=30t+33;$$
$$H_a(17)=Nextprime(30t+34),\enskip H_b(17)=Nextprime(30t+33).$$
So $H_a(n)=H_b(n),\enskip n\geq17,$ and and, for $n>=2,$ $\max(H_a(n)-H_b(n)=6$
which holds for $n=11,12.\blacksquare$

\section{Proof of Theorem 1 (necessity) and of Theorem 2}

1) Let $c(n)=p\equiv 7 \pmod {10}.$ By Section 2, if we have five consecutive
primes of the form
\begin{equation}\label{5}
30t+17=p,30t+19,30t+23,30t+29,30t+31,\enskip t\geq0,
\end{equation}
then, for all $n\geq2,$ \enskip $(H_a(n)-H_b(n))\leq6.$
 Let us show that, if the condition (\ref{5}) does not hold, then $(H_a(n)-H_b(n))>6.$
Denote by $m$ the smallest $n$ when $(H_a(n)-H_b(n))>6.$\newline

1a) Let $30t+17, 30t+19$ be prime, but $30t+23$ be composite. Then we have
$$H_a(2)=30t+19,\enskip H_b(2)=30t+17;$$
$$H_a(3)=Nextprime(30t+23)\geq30t+29,\enskip H_b(3)=30t+19.$$
So $(H_a(3)-H_b(3))\geq10,$ i.e. $m=3.$ \newline

1b) Let $30t+17, 30t+19, 30t+23$ be prime, but $30t+29$ be composite. We have
$$H_a(2)=30t+19,\enskip H_b(2)=30t+17;$$
$$H_a(3)=30t+23,\enskip H_b(3)=30t+19;$$
$$H_a(4)=30t+24,\enskip H_b(4)=30t+20;$$
$$H_a(5)=Nextprime(30t+29)\geq30t+31,\enskip H_b(5)=30t+23.$$
So $(H_a(5)-H_b(5))\geq8,$ i.e. $m=5.$
\newpage
1c) Let $30t+17, 30t+19, 30t+23, 30t+29$ be prime, but $30t+31$ be composite. We have
$$H_a(4)=30t+24,\enskip H_b(4)=30t+20;$$
$$H_a(5)=30t+29,\enskip H_b(5)=30t+23;$$
$$H_a(6)=30t+30,\enskip H_b(6)=30t+24;$$
$$H_a(7)=Nextprime(30t+31)\geq30t+37,\enskip H_b(7)=30t+29.$$
 So $(H_a(7)-H_b(7))\geq8,$ i.e. $m=7.$
  This proves the case $c(n)\equiv 7 \pmod {10}$ with $m=3,5,\enskip or \enskip 7$\newline

2) Now let $c(n)=p\equiv 9 \pmod {10}.$ By Section 2, if we have five consecutive
primes of the form
\begin{equation}\label{6}
30t+29=p,30t+31,30t+37,30t+41,30t+43,\enskip t\geq0,
\end{equation}
then, for all $n\geq2,$ \enskip $(H_a(n)-H_b(n))\leq6.$
 Let us show that, if the condition (\ref{6}) does not hold, then $(H_a(n)-H_b(n))>6$ and $3\leq m\leq17.$
\newline

2a) Let in (\ref{6}) $30t+37$ is composite. Then we have

$$H_a(2)=30t+31,\enskip H_b(2)=30t+29;$$
$$H_a(3)=Nextprime(30t+37)\geq30t+41,\enskip H_b(3)=30t+31.$$
So $(H_a(3)-H_b(3))\geq10,$ i.e. $m=3.$ \newline

2b) Let $30t+37$ be prime, but $30t+41$ be composite. \newline

2ba) Let $30t+37$ be prime, but $30t+41, 30+49$ be composite.\newline
We have

$$H_a(2)=30t+31,\enskip H_b(2)=30t+29;$$
$$H_a(3)=30t+37,\enskip H_b(3)=30t+31;$$
$$H_a(4)=30t+38,\enskip H_b(4)=30t+32;$$
$$H_a(5)=Nextprime(30t+43)\geq30t+47,\enskip H_b(5)=30t+37,$$ if $30t+43$ is composite, $m=5.$
$$H_a(5)=30t+43,\enskip H_b(5)=30t+37,$$ if $30t+43$ is prime;
$$H_a(6)=30t+44,\enskip H_b(6)=30t+38;$$
$$H_a(7)=Nextpime(30t+47)\geq30t+53,\enskip H_b(7)=30t+43,$$
\newpage
 if $30+47$ is composite, $m=7.$
$$H_a(7)=30t+47,\enskip H_b(7)=30t+43,$$ if $30t+47$ is prime;
$$H_a(8)=30t+48,\enskip H_b(8)=30t+44;$$
$$H_a(9)=30t+49,\enskip H_b(9)=30t+45;$$
$$H_a(10)=30t+50,\enskip H_b(10)=30t+46;$$
$$H_a(11)=Nextprime(30t+53)\geq30t+59,\enskip H_b(11)=30t+47;$$ if $30t+53$ is composite, $m=11.$
$$H_a(11)=30t+53,\enskip H_b(11)=30t+47,$$ if $30t+53$ is prime;
$$H_a(12)=30t+54,\enskip H_b(12)=30t+48;$$
$$H_a(13)=Nextprime(30t+54)\geq30t+61,\enskip H_b(13)=30t+53,$$ if $30+59$ is composite, $m=13.$
$$H_a(13)=30t+59,\enskip H_b(13)=30t+53,$$ if $30t+59$ is prime;
$$H_a(14)=30t+60,\enskip H_b(14)=30t+54;$$
$$H_a(15)=30t+62,\enskip H_b(15)=30t+55,$$ if $30t+61$ is prime, $m=15.$
$$H_a(15)=30t+61,\enskip H_b(15)=30t+55,$$ if $30t+61$ is composite;
$$H_a(16)=30t+62,\enskip H_b(16)=30t+56;$$
$$H_a(17)=Nextprime(30t+62)\geq30t+67,\enskip H_b(17)=30t+59,$$ and $m=17,$
which completes proof of the case 2ba)\newline

2bb) Let $30t+37, 30t+49$ be prime, but $30t+41,$ $30t+47$ be composite.\newline
Now the rows $H_a(2),...,H_a(6)$ are as in 2ba). But here we have
 $30t+43$ and $30t+49$ are prime. So,

$$H_a(7)=30t+49,\enskip H_b(7)=30t+43,$$
$$H_a(8)=30t+50,\enskip H_b(8)=30t+44;$$
$$H_a(9)=30t+51,\enskip H_b(9)=30t+46;$$
$$H_a(10)=30t+52,\enskip H_b(10)=30t+47;$$
 \newpage
$$H_a(11)=Nextprime(30t+53)\geq30t+59,\enskip H_b(11)=30t+49,$$ if $30t+53$ is composite, $m=11.$
$$H_a(11)=30t+53,\enskip H_b(11)=30t+49,$$
if $30t+53$ is prime;
$$H_a(12)=30t+54,\enskip H_b(12)=30t+50;$$
$$H_a(13)=Nextprime(30t+59)\geq30t+61,\enskip H_b(13)=30t+53,$$ if $30t+59$ is composite, $m=13.$
$$H_a(13)=30t+59,\enskip H_b(13)=30t+53,$$ if $30t+59$ is prime;
$$H_a(14)=30t+60,\enskip H_b(14)=30t+54;$$
$$H_a(15)=30t+62,\enskip H_b(15)=30t+55,$$ if $30t+61$ is prime, $m=15.$
$$H_a(15)=30t+61,\enskip H_b(15)=30t+55,$$ if $30t+61$ is composite;
$$H_a(16)=30t+62,\enskip H_b(16)=30t+56;$$
$$H_a(17)=Nextprime(30t+62)\geq30t+67,\enskip H_b(17)=30t+59$$ and $m=17,$
which completes proof of the case 2bb)\newline

2bc) Let $30t+37, 30t+47, 30t+49$ be prime, but $30t+41,$ be composite.\newline
Now again the rows $H_a(2),...,H_a(6)$ are as in 2ba). But here we have
 $30t+43,$ $30t+47$ and $30t+49$ are prime. So,
$$H_a(7)=30t+47,\enskip H_b(7)=30t+43,$$
$$H_a(8)=30t+48,\enskip H_b(8)=30t+44;$$
$$H_a(9)=30t+50,\enskip H_b(9)=30t+45;$$
$$H_a(10)=30t+51,\enskip H_b(10)=30t+46;$$
$$H_a(11)=Nextprime(30t+53)\geq30t+59,\enskip H_b(11)=30t+47,$$
if $30t+53$ is composite, $m=11.$
$$H_a(11)=30t+53,\enskip H_b(11)=30t+47,$$
if $30t+53$  is prime;
$$H_a(12)=30t+54,\enskip H_b(12)=30t+48;$$
$$H_a(13)\geq30t+59, \enskip H_b(12)=30t+49$$
\newpage
and $m=13.$\newline
2c) Let $30t+37, 30t+41$ be prime, but $30t+43$ be composite.\newline
We have
$$H_a(2)=30t+31,\enskip H_b(2)=30t+29;$$
$$H_a(3)=30t+37,\enskip H_b(3)=30t+31;$$
$$H_a(4)=30t+38,\enskip H_b(4)=30t+32;$$
$$H_a(5)=30t+41,\enskip H_b(5)=30t+37$$
$$H_a(6)=30t+42,\enskip H_b(6)=30t+38;$$
$$H_a(7)=Nextprime(30t+47)\geq30t+49,\enskip H_b(7)=30t+41,$$ if $30t+47$ is composite, $m=7.$
$$H_a(7)=30t+47,\enskip H_b(7)=30t+41,$$ if $30t+47$ is prime;
$$H_a(8)=30t+48,\enskip H_b(8)=30t+42;$$
$$H_a(9)=30t+50,\enskip H_b(9)=30t+43,$$
if $30t+49$ is prime, $m=9.$
$$H_a(9)=30t+49,\enskip H_b(9)=30t+43,$$ if $30t+49$ is composite;
$$H_a(10)=30t+50,\enskip H_b(10)=30t+44;$$
$$H_a(11)=Nextprime(30t+53)\geq30t+59,\enskip H_b(11)=30t+47,$$ if $30t+53$ is composite, $m=11.$
$$H_a(11)=30t+53,\enskip H_b(11)=30t+47,$$ if $30t+53$ is prime.\newline
The two last rows coincide with the corresponding rows of 2ba). Therefore, 2c) has the same end as 2ba) with $m=13, 15$ and $m=17.$\newline
This proves the case 2c) and together with 2a), 2ba), 2bb) and 2bc) completes
the proof of the case $c(n)\equiv 9 \pmod {10}$ with $m=3,5,7,9,11,13,15 $ 
or \enskip17.\newline

3) Finally let $c(n)=p\equiv 1 \pmod {10}.$ By Section 2, if we have seven consecutive
primes of the form
\begin{equation}\label{7}
30t+11=p,30t+13,30t+17,30t+19,30t+23,30t+29,30t+31,\enskip t\geq0,
\end{equation}
then, for all $n\geq2,$ \enskip $(H_a(n)-H_b(n))\leq6.$
 Again let us show that, if the condition (\ref{7}) does not hold, then $(H_a(n)-H_b(n))>6$ and
 $3\leq m\leq13.$
 \newline
3a) Let $30t+11, 30t+13$ be prime, but $30t+17$ be composite.\newpage
 We have
 $$H_a(2)=30t+13,\enskip H_b(2)=30t+11;$$
 $$H_a(3)=Nextprime(30t+19)\geq30t+23,\enskip H_b(3)=30t+13,$$ if $30t+19$ is composite, $m=3.$
 $$H_a(3)=30t+19,\enskip H_b(3)=30t+13,$$ if $30t+19$ is prime;
 $$H_a(4)=30t+20,\enskip H_b(4)=30t+14;$$
 $$H_a(5)=Nextprime(30t+23)\geq30t+29,\enskip H_b(5)=30t+19,$$ if $30t+23$ is composite, $m=5.$
 $$H_a(5)=30t+23,\enskip H_b(5)=30t+19,$$ if $30t+23$ is prime.
 $$H_a(6)=30t+24,\enskip H_b(6)=30t+20;$$
 $$H_a(7)=Nextprime(30t+29)\geq30t+31,\enskip H_b(7)=30t+23,$$
 if $30t+29$ is composite, $m=7.$
 $$H_a(7)=30t+29,\enskip H_b(7)=30t+23,$$ if $30+29$ is prime;
 $$H_a(8)=30t+30,\enskip H_b(8)=30t+24;$$
 $$H_a(9)=30t+32,\enskip H_b(9)=30t+25,$$ if $30t+31$ is prime and $m=9.$
 $$H_a(9)=30t+31,\enskip H_b(9)=30t+25,$$ if $30t+31$ is composite;
 $$H_a(10)=30t+32,\enskip H_b(10)=30t+26;$$
 $$H_a(11)=Nextprime(30t+32)\geq30t+37,\enskip H_b(11)=30t+29$$
 and $m=11.$ This proves 3a).\newline

3b) Let $30t+11, 30t+13, 30t+17$ be prime, but $30t+19$ be composite.\newline
 We have
 $$H_a(2)=30t+13,\enskip H_b(2)=30t+11;$$
 $$H_a(3)=30t+17,\enskip H_b(3)=30t+13;$$
 $$H_a(4)=30t+18,\enskip H_b(4)=30t+14;$$
 $$H_a(5)=Nextprime(30t+23)\geq30t+29,\enskip H_b(5)=30t+17,$$ if $30t+23$ be composite, $m=5.$
 \newpage
 $$H_a(5)=30t+23,\enskip H_b(5)=30t+17,$$ if $30t+23$ is prime;
 $$H_a(6)=30t+24,\enskip H_b(6)=30t+18;$$
 $$H_a(7)=Next(30t+29)\geq31,\enskip H_b(7)=30t+23,$$ if $30t+29$ be composite, $m=7.$
 $$H_a(7)=30t+29,\enskip H_b(7)=30t+23,$$ if $30t+29$ be prime.
 Since the last two rows coincide with the corresponding rows $H_a(7)'s$ of the
 previous point 3a), then 3b) has the same end as 3a) with $m=9$ and $m=11.$
  This proves 3b).\newline

3c) Let $30t+11, 30t+13, 30t+17, 30t+19$ be prime, but $30t+23$ be composite.\newline
 We have
$$H_a(2)=30t+13,\enskip H_b(2)=30t+11;$$
 $$H_a(3)=30t+17,\enskip H_b(3)=30t+13;$$
 $$H_a(4)=30t+18,\enskip H_b(4)=30t+14;$$
 $$H_a(5)=30t+19,\enskip H_b(5)=30t+17;$$
 $$H_a(6)=30t+20,\enskip H_b(6)=30t+18;$$
 $$H_a(7)=Nextprime(30t+23)\geq30t+29,\enskip H_b(7)=30t+19$$ and $m=7.$
 This proves 3c).\newline

3d) Let $30t+11, 30t+13, 30t+17, 30t+19, 30t+23$ be prime, but $30t+29$ be composite.\newline
We have
$$H_a(6)=30t+20,\enskip H_b(6)=30t+18;$$
$$H_a(7)=30t+23,\enskip H_b(7)=30t+19;$$
$$H_a(8)=30t+24,\enskip H_b(8)=30t+20;$$
$$H_a(9)=30t+25,\enskip H_b(9)=30t+21;$$
$$H_a(10)=30t+26,\enskip H_b(10)=30t+22;$$
$$H_a(11)=Nextprime(30t+29)\geq30t+31,\enskip H_b(11)=30t+23$$ and $m=11.$
This proves 3d). \newline
Finally,\newline
3e) Let $30t+11, 30t+13, 30t+17, 30t+19, 30t+23, 30t+29$ be prime, but $30t+31$ be composite.\newpage
We have
$$H_a(10)=30t+26,\enskip H_b(10)=30t+22;$$
$$H_a(11)=30t+29,\enskip H_b(11)=30t+23;$$
$$H_a(12)=30t+30,\enskip H_b(12)=30t+24;$$
$$H_a(13)=Nextprime(30t+31)\geq30t+37,\enskip H_b(13)=30t+29$$ and $m=13.$
This completes the proof of the case 3e) and the proof of the case $c(n)\equiv 1 \pmod {10}$ with $m=3,5,7,9,11,\enskip or\enskip13.$ The proved cases 1),2) and 3) complete the proof
of both Theorem \ref{t1} and \ref{2}. $\blacksquare$
\begin{corollary}\label{c1}
 In sequence $(\ref{2})$ there is only one term equal to $4$; all other terms are
greater than or equal to $6.$
\end{corollary}

Also from 2ba), 2bb) and 2c) we deduce the following.

\begin{corollary}\label{c2}
In sequence $(\ref{4})$  $m=17$ appears if and only if we have seven primes:
either
$$p=30t+29,p+2,p+8,p+12,p+18,p+24,p+30$$
or
$$p,p+2,p+8,p+14,p+20,p+24,p+30$$
or
$$p,p+2,p+8,p+14,p+18,p+24,p+30,$$
but $p+32$ is composite.
\end{corollary}
\begin{corollary}\label{c3}
In sequence $(\ref{4})$  $m=15$ appears if and only if we have one of the constellations of eight primes: either
$$p=30t+29,p+2,p+8,p+12,p+18,p+24,p+30,p+32$$
or
$$p,p+2,p+8,p+14,p+20,p+24,p+30,p+32$$
or
$$p,p+2,p+8,p+14,p+18,p+24,p+30,p+32.$$
\end{corollary}
So $m=15$ occurs asymptotically less than $m=17.$
Analogous statements one can obtain for other positive values of $m$ among which
 $m=3$ appears asymptotically most frequently.

\end{document}